\documentclass[12pt, twoside]{article}
\usepackage{amsmath,amsthm,amssymb}
\usepackage{times}
\usepackage{enumerate}
\usepackage{color}
\usepackage{srcltx}

\makeatletter
\def\author#1{\gdef\autrun{\def\and{\unskip, }#1}\gdef\@author{#1}}

\makeatother

\def\subjclass#1{{\renewcommand{\thefootnote}{}%
\footnote{\emph{Mathematics Subject Classification (2010):} #1}}}
\def\keywords#1{\par\medskip
\noindent\textbf{Keywords.} #1}

\def\CX{{\mathbb C}}

\def\NX{{\mathbb N}}
\def\ZX{{\mathbb Z}}

\def\GL{{\rm GL}}

\def\gl{{\rm gl}}

\def\QED{\hbox{\hskip 1pt \vrule width4pt height 6pt depth 1.5pt \hskip 1pt}}

\def\calM{{\cal M}}

\def\calA{{\cal A}}

\newcommand{\ie}{{\it i.e.}}
\newcommand{\diag}{ {\rm diag} }

\renewcommand{\Re}{\re}
\renewcommand{\Im}{\mbox{\rm Im\,}}

\DeclareMathOperator{\re}{Re}

\newtheorem{theorem}{Theorem}
\newtheorem{prop}[theorem]{Proposition}

\newtheorem{lemma}[theorem]{Lemma}

\def\QED{\hbox{\hskip 1pt \vrule width4pt height 6pt depth 1.5pt \hskip 1pt}}
\newenvironment{prff}[1]{\trivlist
\item[\hskip \labelsep{\bf #1.\hspace*{.3em}}]}{~\hspace{\fill}~$\square$\endtrivlist}
\newenvironment{prf}{
\begin{prff}{Proof}}{
\end{prff}}
 \def\square{\QED}

\begin{document}

\title{Mahler equations and rationality}
\author{Reinhard Sch\"afke\footnote{ 
Institut de Recherche Math\'ematique Avanc\'ee,
Universit\'e de  Strasbourg et C.N.R.S.,
7, rue Ren\'e Descartes,
67084 Strasbourg Cedex, France, {\tt schaefke@unistra.fr}.},   
Michael F. Singer\footnote{Department of Mathematics, North Carolina State University, Box 8205, Raleigh, NC 27695, USA, {\tt singer@ncsu.edu}.  }}

\date{March 21, 2017}
\maketitle

\begin{abstract}{
We give another proof of a result of Adamczewski and Bell \cite{AuB16} concerning Mahler equations:
A formal power series 
satisfying a $p$- and a $q$-Mahler equation over $\CX(x)$ with multiplicatively independent
positive integers $p$ and $q$ is a rational function.
The proof presented here is self-contained and 
essentially a compilation of proofs contained in a recent preprint \cite{MSRS16} of the authors.
}\end{abstract}
We consider two Mahler operators, \ie\ two endomorphisms $\sigma_j$, $j=1,2$,
on the field $K=\CX[[x]][x^{-1}]$ of formal Laurent series with complex coefficients
defined by $\sigma_1(f(x))=f(x^{p})$, $\sigma_2(f(x))=f(x^{q})$ for any $f(x)\in K$ where $p$ and $q$
are positive integers. Observe that $\sigma_1$ and $\sigma_2$ commute, 
\ie\ $\sigma_1\sigma_2=\sigma_2\sigma_1$.
We consider the field $\CX(x)$ of rational functions
with complex coefficients as a subfield of $K$, the inclusion given by
the expansion in a Laurent series at the origin.
We want to prove the following theorem
\begin{theorem} Assume that $p$ and $q$ are multiplicatively independent,
\ie\ there are no nonzero integers $n_j$ such that $q^{n_2}=p^{n_1}$.
Suppose that the formal series $f(x)\in K$ satisfies
a system of two Mahler equations
\begin{equation}\label{sys2s}S_j(f(x))  =  \sigma_j^{m_j}(f(x)) + b_{j,m_j-1}(x) 
\sigma_j^{m_j-1}(f(x)) + \ldots + b_{j,0}(x)f(x) =  0,\ j=1,2 
\end{equation}
with $b_{j,i}(x) \in \CX(x)$.

Then $f(x)$ is rational.\end{theorem}

\noindent {\bf Remark:} 1.\ This Theorem was recently proved  by Adamczewski and Bell in \cite{AuB16}. 
Their tools include a local-global principle to reduce the problem to a similar problem 
over finite fields, Chebotarev's Density Theorem, Cobham's Theorem and  some asymptotics - all very different from the techniques used in the present work.\\
\noindent 2.\ \cite{AuB16} also provides background information about Mahler equations, in particular
historical, many references to the literature and explains the relation to Cobham's theorem in the theory
of finite state machines. The fact that the generating functions of $p$-regular (and thus of $p$-automatic sequences)
satisfy $p$-Mahler equations is shown in \cite{Becker94}.\\ 
\noindent 3.\ The subsequent proof is essentially a compilation of work contained in our
recent preprint \cite{MSRS16}, see Corollary 15, part 3, and Proposition 19. The preprint
presents a unified reduction theory 
of consistent pairs of first order systems of linear differential, difference, 
$q$-difference or Mahler equations like the one of Proposition \ref{2mprop} below and
uses it to deduce numerous statements on common solutions of two scalar linear differential, difference, 
$q$-difference or Mahler equations.

\begin{prf}
We may assume without loss of generality that $b_{j,0}(x)\neq0$, $j=1,2$. This follows from
\begin{lemma}\label{linind}
Consider $w_1(x),...,w_\ell(x)\in K$ and a positive integer $m$. 
Then these series are $\CX(x)$-linearly dependent
if and only if $w_1(x^m),...,w_\ell(x^m)$ are.\end{lemma}
\begin{prf} We only prove the nontrivial implication.
Suppose  that $w_1(x^m),...,w_\ell(x^m)$ are $\CX(x)$-linearly dependent, which means
that there exist $a_k\in\CX(x)$, $k=1,...,\ell$, not all
zero, such that 
$$a_1(x)w_1(x^{m})+...+a_\ell(x)w_\ell(x^{m})=0.$$
Now we can uniquely write $a_k(x)=\sum_{j=0}^{m-1} x^j c_k^{j}(x^{m})$, $k=1,...,\ell$,
with rational functions $c_k^j(x)$. Expanding the terms in the above equation in Laurent series
we obtain the equations
$$c_1^j(x^{m})w_1(x^{m})+...+c_\ell^j(x^{m})w_\ell(x^{m})=0, j=0,...,m-1,$$
and hence
$$c_1^j(x)w_1(x)+...+c_\ell^j(x)w_\ell(x)=0, j=0,...,m-1.$$
Since at least one of them must be nontrivial we obtain
the linear dependence of the $w_j$, $j=1,...,\ell$.~\end{prf}

Consider now the $\CX(x)$-subspace $W$ of $K$ generated by $\sigma_1^m\sigma_2^r(f)$, 
$m=0,...,m_1-1,\ r=0,...m_2-1$. By (\ref{sys2s}), $W$ is invariant under $\sigma_1$ and $\sigma_2$; 
here the fact that the $\sigma_j$ commute is used.

Let $g_1,...,g_n$ be {a}   $\CX(x)$-basis of $W$ with $g_1=f$ and let $g=(g_1,...,g_n)^T$. Then we have
that 
\begin{equation}\label{sysW} \sigma_1(g)=A(x)g,\ \sigma_2(g)=B(x)g,\end{equation}
with $A,B\in \gl_n(\CX(x))$.
By Lemma \ref{linind}, we actually have $A,B\in \GL_n(\CX(x))$ because the {components} of $\sigma_j(g)$ 
form again a basis of $W$.

Additionally, the coefficient matrices of (\ref{sysW}) satisfy a certain {\em consistency condition}. 
Indeed, we have
$$0=\sigma_1(\sigma_2(g))-\sigma_2(\sigma_1(g))=(\sigma_1(B)A-\sigma_2(A)B)g$$
and as the {components} of $g$ form a basis we obtain 
\begin{equation}\label{2mcons}A(x^q)B(x)=B(x^p)A(x).\end{equation}

Our statement  then follows from
\begin{prop}\label{2mprop}Consider a system
\begin{equation}\label{2msys}y(x^p)=A(x)y(x),\ \ \ y(x^q)=B(x)y(x)\end{equation}
with multiplicatively independent positive integers $p$ and $q$ and
$A(x),B(x)\in\GL_n(\CX(x))$ satisfying the consistency condition
(\ref{2mcons}).
Suppose that $g(x)\in (\CX[[x]][x^{-1}])^n$ is a formal vectorial solution. Then
$g(x)\in \CX(x)^n$.\end{prop}

Observe that we must actually have $n=1$ in the proof of the Theorem because
the components of $g(x)$ are $\CX(x)$-linearly independent.\end{prf}

The proof of Proposition~\ref{2mprop} proceeds in three steps.  We first prove that $g(x)$ converges in a neighborhood of $0$.  In the second  step (the heart of the proof) we show that $g(x)$ can be extended analytically to a meromorphic function on $\CX$ with only finitely many poles.  Finally we prove that $g(x)$ has polynomial growth as $|x| \rightarrow \infty$ and therefore must be in $\CX(x)^n$. We begin with the first step.

\begin{lemma}\label{conv} The series $g(x)$ is convergent in a neighborhood of $0$.\end{lemma}

\begin{prf} This is a special case of \cite{Bez94}, Theorem 1-2, and 
could also be deduced from \cite{Ran92}, section 4. 
For the convenience of the reader, we provide a short proof. 
To do that, we truncate $g(x)$ at a sufficiently high power of $x$ to obtain
$h(x)\in (\CX[x][x^{-1}])^n$ and introduce $r(x)=h(x)-A(x)^{-1}h(x^p)$
and $\tilde g(x)=g(x)-h(x)$. Then we have
\begin{equation}\label{2m-gr}
\tilde g(x)=A(x)^{-1}\tilde g(x^p) - r(x).\end{equation}
We denote the valuation of $A(x)^{-1}$ at the origin by $s\in\ZX$
and introduce $\tilde A(x)=x^{-s}A(x)^{-1}$ which is holomorphic at the origin.

First choose $M\in\NX$ such that $pM+s>M$ and $h(x)$ such that $g(x)-h(x)$ has at least
valuation $M$. Then by (\ref{2m-gr}), $r(x)$ also has at least valuation $M$.
Now consider $R>0$ such that $\tilde A(x)$ is holomorphic and bounded on $D(0,R)$. 
Then consider for positive $\rho<\min(R,1)$ the vector space $E_\rho$ of all series  
$F(x)=\sum_{m=M}^\infty F_mx^m$ such that $\sum_{m=M}^\infty|F_m|{\rho}^m$ 
converges and define the norm
$|F(x)|_\rho$ as this sum. Then $E_\rho$ equipped with $|\ |_\rho$ is a Banach space
and the existence of a unique solution of (\ref{2m-gr}) in $E_\rho$
for sufficiently small $\rho>0$ follows from the Banach fixed-point theorem using
that $|x^s F(x^p)|_\rho\leq \rho^{Mp+s-M}|F(x)|_\rho$ for $F(x)\in E_\rho$.  Since any solution 
$y(x) \in x^M\CX[[x]]$ of $y(x) = A(x)^{-1} y(x^p)$ must be zero, we have that $\tilde{g}(x)$ coincides with the solution in $E_\rho$.
This proves the convergence of $\tilde g(x)$ and hence of $g(x)$.
\end{prf} 


We now turn to the task of showing that $g(x)$ can be extended to a meromorphic function on $\CX$.
By (\ref{2msys}), rewritten $g(x)=A(x)^{-1}g(x^p)$, the function $g$ can only be extended
analytically to a meromorphic function on the unit disk. As we want to extend it {beyond}
the unit disk, we use the change of variables $x=e^t$, $u(t)=y(e^t)$ and obtain a system of
$q$-difference equations
\begin{equation}\label{2msys2q}u(pt)=\bar A(t)u(t),\ \ \ u(qt)=\bar B(t)u(t)\end{equation}
with $\bar A(t)=A(e^t)$, $\bar B(t)=B(e^t)$. It satisfies the consistency condition
\begin{equation}\label{2mconsi2q}\bar A(qt)\bar B(t)=\bar B(pt)\bar A(t).\end{equation}
Observe that $\bar A(t),\bar B(t)$ are not rational in $t$, but rational 
in $e^t$.

The heart of the proof of Proposition~\ref{2mprop} lies in understanding the behavior of solutions of (\ref{2msys2q}).  We do this by first showing in Lemma~\ref{const} that there is a formal gauge transformation $u=G\,v$, $G\in\GL_n(\CX\{t\}[t^{-1}])$,  such that $v$ satisfies a system with constant coefficients.  We then show in Lemma~\ref{2mcont} that the transformation matrix $G(t)$ and its inverse can be continued analytically  to meromorphic functions on the $t$-plane. The ``quotient'' function $d(t) = G(t)^{-1}g(e^t)$ then satisfies a system with 
constant coefficients which can be solved explicitly. In this way, we show in Lemma~\ref{d(t)} that $d(t)$  can  be extended analytically to an 
entire function on the Riemann surface of $\log(t)$.  Using these three lemmas, we show in Lemma~\ref{meromg} that $g(x)$ can be continued analytically to a meromorphic function on the $x$-plane.

\begin{lemma}\label{const} There exists a convergent gauge transformation $u=G(t)\,v$, 
$G(t)\in\GL_n(\CX\{t\}[t^{-1}])$, such that $v$ satisfies
\begin{equation}\label{consteq}v(pt) = A_1 v(t) , \ \ \  v(qt) = B_1v(t)
\end{equation}
where $A_1, B_1 \in \GL_n(\CX)$ commute.
\end{lemma}

\begin{prf} Concerning the behavior at $t=0$, it is known that there exists a formal gauge transformation
$u=G\,z$, $G\in\GL_n(\CX[[t^{1/s}]][t^{-1/s}])$, {$s\in\NX^*$,} that reduces $u(pt)=\bar A(t)u(t)$
to a system $z(pt)=t^DA_1z(t)$, where $D$ is a diagonal matrix with entries in $\frac1s\ZX$ and $A_1\in\GL_n(\CX)$ 
such that any eigenvalue $\lambda$ of $A_1$ satisfies $1\leq|\lambda|<|p|^{1/s}$, moreover $D$ and $A_1$ commute.
If we write $D=\diag(d_1I_1,...,d_rI_r)$ with distinct $d_j$ and $I_j$ identity matrices of an appropriate size,
then $A_1=\diag(A_1^{1},...,A_1^r)$ with diagonal blocks $A_1^{j}$ of corresponding size. 
$D$ and $A_1$ are essentially unique, \ie\ except for a permutation of the diagonal blocks and passage from
some $A_1^j$ to a conjugate matrix. If $D$ happens to be 0, then $s$ can
be chosen to be 1 and $G$ is convergent
(see \cite{PuSi}, ch.\ 12, \cite{Ad31}, \cite{Car12}).

Now by the consistency condition (\ref{2mconsi2q}), the gauge transformation $v=B(t)u$ transforms
$u(pt)=\bar A(t)u(t)$ to
$v(pt)=\bar A(qt)v(t)$. The gauge transformation $v=G(qt)w$ then transforms this system to
$w(pt)=(qt)^DA_1w(t)$. Now $(qt)^DA_1=t^{D}\,q^{D}A_1$ and there is a diagonal matrix $F$
with entries in $\frac1s\ZX$ commuting with $D$ and $A_1$ such that the gauge transformation 
$ w=t^{F}\tilde w$ reduces
the latter system to $\tilde w(pt)=t^D\tilde A_1\tilde w(t)$, where $\tilde A_1=p^{-F}q^DA_1$ has again eigenvalues
with modulus in $[1,|p|^{1/s}[$. Now we write $\tilde A_1=\diag(\tilde A_1^{1},...,\tilde A_1^r)$ and 
fix some $j\in\{1,...,r\}$. If $a^{j}_1,...,a^{j}_{r_j}$ are the eigenvalues of $A_1^j$ then
$p^{-f_j}q^{d_j}a^j_\ell$, $\ell=1,...,r_j$, are those of $\tilde A_1^j$. By the uniqueness of the reduced form, 
the mapping $t\mapsto p^{-f_j}q^{d_j}t$ induces a permutation of the eigenvalues of $A_1^j$. 
If we apply it several times, if necessary, we obtain the existence of some $\ell\in\{1,...,r_j\}$ and of some positive integer
$k$ such that $p^{-kf_j}q^{kd_j}a^j_\ell=a_\ell^j$. Due to our condition on $p$ and $q$ this is only possible
if $d_j=0$. Thus we have proved that $D=0$ and $t=0$ is a so-called {\em regular singular point}
of $u(pt)=\bar A(t)u(t)$.

We {therefore} obtain a matrix $A_1$ with eigenvalues $\lambda$ in the annulus $1\leq|\lambda|<p$
and $G(t)\in\GL_n(\CX\{t\}[t^{-1}])$ such that $u=G(t)v$ reduces the first equation of (\ref{2msys2q})
to $v(pt)=A_1v(t)$. This means 
\begin{equation}\label{2mg1}G(pt)=\bar A(t)G(t)A_1^{-1}\mbox{ for small }t.\end{equation}
Applying the same gauge transformation to the second equation of (\ref{2msys2q})
yields an equation $v(qt)=\tilde{\bar B}(t)v(t)$ with some $\tilde{\bar B}(t)\in\GL_n(C\{t\}[t^{-1}])$.
It satisfies the consistency condition $A_1\tilde{\bar B}(t)=\tilde{\bar B}(pt)A_1$.
Now we expand $\tilde{\bar B}(t)=\sum_{m=m_0}^\infty C_m t^m$.
The coefficients satisfy $A_1C_m=C_m (p^mA_1)$, $m\geq m_0$. As $A_1$ and $p^mA_1$ have no common eigenvalue
unless $m=0$, we obtain that $\tilde{\bar B}(t)=:B_1$ is constant and commutes with $A_1$.
We note the second equation satisfied by $G$  
\begin{equation}\label{2mg1b}G(qt)=\bar B(t)G(t)B_1^{-1}\mbox{ for small }t.\end{equation} \end{prf}

{\begin{lemma}\label{2mcont}
The functions $G(t)^{\pm1}$ 
can be continued analytically to meromorphic functions on $\CX$
and there exists $\delta>0$ such that both  can be continued analytically
to the sectors $\{t\in\CX^*\mid \delta<\arg(\pm t)<2\delta\}$.
\end{lemma}
\begin{prf}
Let $\calM$ be the set of poles of $\bar A(t)^{\pm1}$, \ie\ the set of $t$
such that $e^t$ is a pole of $\bar A(x)$ or $\bar A(x)^{-1}$. Note that $\calM$ is $2 \pi i$-periodic,
has no finite accumulation point and is contained in some vertical strip
$\{t\in\CX\mid -D<\Re t<D\}$.
By (\ref{2mg1}), $G(t)^{\pm1}$ can be continued analytically to $\CX^*\setminus(\calM\cdot p^\NX)$
and thus to meromorphic functions on $\CX$ which we denote by the same name. By construction,
$G(t)^{\pm1}$ are also analytic in some punctured neighborhood of the origin.
By the properties of $\calM$, the infimum of the $|\Re t_1| $ on the set of all $t_1\in\calM$
having nonzero real part is a positive number. As $\calM$ is contained in some vertical strip
there exist sectors $\{t\in\CX^*\mid \delta<\arg(\pm t)<2\delta\}$ disjoint to $\calM$
and hence to $\calM\cdot p^\NX$. Therefore $G(t)^{\pm1}$ can be analytically continued to these sectors
and the lemma is proved.\end{prf}

\begin{lemma} \label{d(t)} The function $d(t)=G(t)^{-1}g(e^t)$ can be continued analytically to the Riemann surface of $\log(t)$. \end{lemma}

\begin{prf} By Lemma \ref{2mcont} and because $g(x)$
is holomorphic in some punctured neighborhood of $x=0$ by Lemma \ref{conv}, $d(t)$ is defined
and holomorphic for some sector
$S=\{t\in\CX\mid |t|>K,\ \pi+\delta<\arg t<\pi+2\delta\}$. By (\ref{2msys}), (\ref{2mg1}),
and (\ref{2mg1b}) it satisfies
\begin{equation}\label{2md}d(pt)=A_1d(t),\ \ d(qt)=B_1d(t)\mbox{ for }t\in S.\end{equation}
To solve (\ref{2md}), consider a matrix 
$L_1$ commuting with $B_1$ such that $p^{L_1}=A_1$. Put $F(t)=t^{-L_1}d(t)$.
Then
\begin{equation}\label{2mf}F(pt)=F(t),\ \ F(qt)=\tilde B_1F(t)\mbox{ for }t\in S\end{equation}
where {$\tilde B_1=B_1q^{-L_1}$}. Thus $H(s)=F(e^s)$ is $\log(p)$-periodic on the half-strip
$B=\{s\in\CX\mid \Re s>\log(K),\ \pi+\delta<\Im s<\pi+2\delta\}$ and can be expanded in a Fourier series.
This implies that
\begin{equation}\label{2mfourier}F(t)=\sum_{\ell=-\infty}^\infty F_\ell\,t^{\frac{2 \pi i}{\log(p)}\ell}\mbox{ for }t\in S.
\end{equation}
The second equation of (\ref{2mf}) yields conditions on the Fourier coefficients
$$F_\ell\exp\left(2\pi i\tfrac{\log(q)}{\log(p)}\ell\right)=\tilde B_1F_\ell\mbox{ for }\ell\in\ZX.$$
Therefore $F_\ell=0$ unless $\exp\left(2\pi i\frac{\log(q)}{\log(p)}\ell\right)$ is an eigenvalue of 
$\tilde B_1$.
Since $p$ and $q$ are multiplicatively independent, the quotient $\frac{\log(q)}{\log(p)}$ is irrational and hence
$\exp\left(2\pi i\frac{\log(q)}{\log(p)}\right)$ is not a root of unity. Therefore all the numbers
$\exp\left(2\pi i\frac{\log(q)}{\log(p)}\ell\right)$, $\ell\in\ZX$ are different and only finitely many of them
can be eigenvalues of $\tilde B_1$. This shows that the Fourier series  (\ref{2mfourier}) has finitely many terms
and thus $F(t)$ can be analytically continued to the whole Riemann surface $\hat \CX$ of $\log(t)$. The same holds
for $d(t)=t^{L_1}F(t)$. \end{prf}

\begin{lemma} \label{meromg} The function $g(x)$ can be continued analytically to a meromorphic function on $\CX$ with finitely many poles.
\end{lemma}

\noindent{\bf Remark:} According to Theorem 4.2 of \cite{Ran92} (see also \cite{BCR13}), it is sufficient to show that $g(x)$ {does not have} the
unit circle as a natural boundary and the rationality of $g(x)$ follows. {We show how it follows naturally, in our context,} that $g(x)$ can be continued 
analytically as a meromorphic function to all of $\CX$ and, as well, that it has only finitely many poles.
The rationality of $g(x)$ then follows as in  \cite{Ran92} and \cite{BCR13} from a growth estimate (Lemma \ref{grow}).

\begin{prf} The function $h(t)=g(e^t)$ is holomorphic for $t$ with large negative real part by Lemma \ref{conv} and $2\pi i$-periodic.
Using Lemma \ref{2mcont} we conclude that $h(t)=G(t)d(t)$ can be
analytically continued to a meromorphic function on $\hat \CX$, in particular the point $t=2\pi i$ is at most a pole of
$h$. By its periodicity, this implies that
$t=0$ also is at most a pole of $h$ and that
it can be continued analytically to a meromorphic function on $\CX$  which we denote by the same name.


Since $h(t)=g(e^t)$ for  $t$ with large negative real part, $h(t)$ is $2\pi i$-periodic for those
values of $t$, hence also its analytic continuation to a meromorphic function on all of $\CX$.
This periodicity allows {one} to define a meromorphic function $\tilde g(x)$ on {$\CX\backslash\{0\}$} by
$\tilde g(e^t)=h(t)$. As $\tilde g(x)=g(x)$ for small $|x|\neq0$ by the construction of $h$, we
have shown that $g(x)$ can be continued analytically to a meromorphic function on $\CX$
which will again be denoted by the same name.

The formula $h(t)=G(t)d(t)$ and Lemma \ref{2mcont} also imply that $h$ is analytic in the sector
$\tilde S=\{t\in\CX^*\mid \delta<\arg t<2\delta\}$. As this sector contains some half strip
$\{t\in\CX\mid \Re t>L,\mu\Re t<\Im t<\mu\Re t+3\pi\}$ 
for some positive $L,\mu$ which has vertical width larger than $2\pi$ and $h$ is $2\pi i$-periodic,
its poles are contained in some vertical strip $\{t\in\CX\mid -L<\Re t<L\}$. This implies that  $g(x)$ has only a finite number of poles.\end{prf}
}

The proof of Proposition~\ref{2mprop} is completed once we  have shown

\begin{lemma}\label{grow}The function  $g(x)$ has polynomial growth as $|x|\to\infty$. \end{lemma}

\begin{prf} This is shown in 
the proof of Theorem 4.2 in \cite{Ran92} (see also \cite{BCR13}). 
For the convenience of the reader, we reproduce it below.

Consider
$r_0{>1}$ such that $g(x)$ and $A(x)$ are holomorphic on the annulus $|x|>r_0/2$. There are
positive numbers $K,M$ such that $|A(x)|\leq K|x|^M$ for $|x|\geq r_0$.
Consider now the annuli
$$\calA_j=\{x\in\CX\mid r_0^{p^j}\leq|x|<r_0^{p^{j+1}}\},\ j=0,1,...$$
covering the annulus $|x|\geq r_0$. Any $x\in\calA_j$ can be written $x=\xi^{p^j}$ with some $\xi\in\calA_0$.
Then we estimate using (\ref{2msys}) and the inequality for $|A(x)|$
$$|g(x)|=|g(\xi^{p^j})|\leq K ^j \left(|\xi|^{p^{j-1}}\cdots|\xi|^p|\xi|\right)^M \max_{r_0\leq|\xi|\leq r_0^p}|g(\xi)|.$$
Hence there is a positive constant $L$ such that
$|g(x)|\leq L\, K^j \,|x|^{\frac M{p-1}}\mbox{ for }x\in\calA_j.$
Assuming $\log(r_0)\geq1$ without loss in generality, we find that $j\leq \log(\log(|x|))/\log(p)$ for $x\in\calA_j$.
Hence there exists $d>0$ such that
$$|g(x)|\leq L\, (\log(|x|))^d\,|x|^{\frac M{p-1}}\mbox{ for }|x|>r_0.$$ \end{prf}

{ \noindent{Acknowledgement.} The authors would like to thank Boris Adamczewski for suggesting an improvement of the 
proof of Lemma \ref{meromg} and for pointing out the article \cite{BCR13} to us.}

\bibliographystyle{plain}
\bibliography{DDrefs}

\end{document}